\documentclass[10pt]{amsart}
\usepackage{a4,a4wide}
\usepackage{leqno}
\usepackage{graphics}
\usepackage{amssymb}
\usepackage{amsmath}
\usepackage{amsthm}
\usepackage{palatino}
\usepackage[T1]{fontenc}

\usepackage{color}

\newtheorem{theorem}{Theorem}[section]

\theoremstyle{definition}

\newtheorem{proposition}[theorem]{Proposition}
\theoremstyle{remark}

\newtheorem{cla}{\textbf{Claim --}}[section]

\numberwithin{equation}{section}

\def\O{{\Omega}}
\def\o{{\omega}}

\def\eps{{\epsilon}}

\def\e{{\mathcal{E}}}

\def\k{{\mathcal{K}}}

\def\L{{\mathcal{L}}}

\def\M{{\mathcal{M}}}

\def\R{{\mathbb{R}}}

\newcommand{\ope}[1]{\e[{#1}]}
\newcommand{\opk}[1]{\k[{#1}]}
\newcommand{\opks}[1]{\k^*[{#1}]}

\newcommand{\lb}[1]{\L_{_{#1}}}
\newcommand{\oplb}[2]{\L_{_{#2}}[{#1}]}
\newcommand{\lbs}[1]{\L^*_{_{#1}}}
\newcommand{\oplbs}[2]{\L^*_{_{#2}}[{#1}]}

\newcommand{\dem}[1]{\vskip 0.2\baselineskip \noindent {\bf{#1}}\vskip 0.2\baselineskip }
\newcommand{\fdem}{\vskip 0.2 pt \qquad \qquad \qquad \qquad \qquad \qquad \qquad \qquad \qquad \qquad \qquad \qquad \qquad \qquad \qquad \qquad  $\square$  }

\begin{document}

\title{Nonlocal refuge model with a partial control}

\author{J{\'e}r{\^o}me Coville\\
}
\address{UR 546 Biostatistique et Processus Spatiaux\\
 INRA, Domaine St Paul Site Agroparc\\
  F-84000 Avignon\\
 France\\}
\email{jerome.coville@avignon.inra.fr}
\thanks{The author is supported by INRA Avignon }

\subjclass[2000]{ }

\date{May 30, 2013.}


\begin{abstract}
In this paper, we analyse the structure of the set of  positive solutions of an heterogeneous nonlocal equation of the form:
$$  \int_{\O} K(x, y)u(y)\,dy -\int_ {\O}K(y, x)u(x)\, dy + a_0u+\lambda a_1(x)u -\beta(x)u^p=0 \quad \text{in}\quad  \times \O$$
where $\O\subset \R^n$ is a bounded open set,  $K\in C(\R^n\times \R^n) $ is nonnegative, $a_i,\beta \in C(\O)$ and $\lambda\in\R$. 
 Such type of equation appears in some studies of  population dynamics where the above solutions  are the stationary states of the dynamic of a spatially structured population evolving in a heterogeneous partially controlled landscape and submitted to a long range dispersal.
 Under some fairly general assumptions on $K,a_i$ and $\beta$ we first establish a necessary and sufficient criterium for the existence of a unique positive solution. Then we analyse the structure of the set of positive solution $(\lambda,u_\lambda)$ with respect to the presence or absence of a refuge zone (i.e $\o$ so that $\beta_{|\o}\equiv 0$).  
\end{abstract}

\maketitle
\section{Introduction}
In this article we are interested in the positive bounded solutions of the nonlinear nonlocal equation 
 \begin{equation}
  \int_{\O} K(x, y)u(y)\,dy -k(x)u + a_0(x)u+\lambda a_1(x)u -\beta(x)u^p=0 \quad \text{in}\quad  \O, \label{refuge-eq}
 \end{equation}
where $\O\subset \R^n$ is a bounded open set, $K\in C(\R^n\times \R^n) $ is non negative, $k(x):=\int_{\O}K(y,x)\,dy$ $\lambda \in \R$, and $a_i,\beta$ are continuous functions. Our aim is to describe the properties of the positive bounded solutions of \eqref{refuge-eq-parab}, in terms of the properties of $K, a_i,\beta$ and $\lambda$.  That is, we look for existence  criteria  of  positive bounded solutions of \eqref{refuge-eq} and we describe some bifurcation diagrams i.e. depending on $a_i$ and $\beta$ we analyse the properties of the curve $(\lambda,u_\lambda)$.

The study of these kind of problems finds its justification in the ecological problematics related to the erosion of Biodiversity. In particular, some recent studies have focused on  a better understanding of the impact of some agricultural practises on non targeted species \cite{Arpaia2010,Birch2005,Hendriksma2012,Perry2012,Perry2010,Pleasants2013}. Such problematic can  be addressed through the analysis of the asymptotic behaviour of the positive solution of a reaction diffusion equation :
\begin{align} 
 &\frac{\partial u(t,x)}{\partial t} = \int_{\O} K(x, y)u(t,y)\,dy -k(x)u(t,x) + a_0(x)u+\lambda a_1(x)u -\beta(x)u^p \quad \text{in}\quad  \R^+\times \O\label{refuge-eq-parab}\\
 &u(0,x) = u_0 (x) \quad \text{ in }\quad \O  \label{refuge-eq-parab-ci}
 \end{align}
where $u$ represents a population density  evolving in a partial controlled heterogeneous . Here the parameter $\lambda$ is  a  control related to the practise and $a_1$ represents the region where the control is exerted.

In the literature the characterisation of the positive bounded  solutions has been extensively studied for the elliptic equations 
\begin{align} 
& \ope{u} + a_0 u+\lambda a_1(x)u -\beta(x)u^p=0 \quad \text{in}\quad  \O, \label{refuge-eq-rd}\\
 &u(x)=0, \quad \text{ in }\quad\partial\O.  \label{refuge-eq-rd-bc}
 \end{align}
   where $\ope{u}:=a_{ij}(x)\partial_{ij}u +b_i(x)\partial_i u + c(x)$ is uniform elliptic  \cite{Berestycki1981,Berestycki2005a,Berestycki2007,Berestycki2006,Cantrell1991a,Cantrell1989,Cantrell1991,Garcia-Melian1998,Garcia-Melian2001a,Ouyang1992}. 
Nowadays, the structure of the positive bounded solutions $u_\lambda$  to \eqref{refuge-eq-rd}--\eqref{refuge-eq-rd-bc} is well understood. More precisely, a positive bounded  solution $u$ to  \eqref{refuge-eq-rd}-- \eqref{refuge-eq-rd-bc} exists \textbf{if and only if}
   \begin{equation}
   \mu_1(\e +a_0+\lambda a_1,\O)<0<\mu_1(\e +a_0+\lambda a_1,\o), \label{refuge-eq-rd-cara}
   \end{equation}
   where $\o$ denotes the refuge zone, i.e. $\o:=\{x\in\O\,|\, \beta(x)=0\}$ and $\mu_1(\O)$ denotes  the  first eigenvalue  of the spectral problem $\ope{\phi}+ a_0 \phi+\lambda a_1(x)\phi +\mu \phi=0$, $\phi=0\; \text{ on }\; \partial \O$.  
   Depending on the properties of $\beta$ and $a_1$ a  description of the curves $(\lambda,u_\lambda)$ can be found in \cite{Fraile1996,Garcia-Melian1998,Garcia-Melian2001a,Ouyang1992}.  
   
For nonlocal equations such as \eqref{refuge-eq},  less is known and the analysis of the existence, uniqueness and the bifurcation diagram have been only studied in particular situations \cite{Bates2007,Coville2010,Coville2008,Coville2012,Garcia-Melian2009a,Hutson2003,Li2012,Shen2012}.  A large part of the literature is devoted to the existence of positive solution to \eqref{refuge-eq} in situations where no refuge zone exists and for a fixed $\lambda$ \cite{Bates2007,Coville2010,Coville2008,Coville2012,Hutson2003,Shen2012}. To our knowledge \cite{Garcia-Melian2009a} is the first paper which considers a nonlocal logistic equation  with a refuge zone and  analyses the curves $(\lambda,u_\lambda)$.  More precisely, the authors investigate the  existence, uniqueness of a positive bounded solution of 
\begin{align}
&J\star u -u +\lambda u -\beta(x) u^p=0 \quad \text{ in }\quad \O, \\
&u\equiv 0 \quad \text{ in }\quad \R^n\setminus \bar\O,
\end{align}
where $J$ is a symmetric density of probability. They prove that a  positive solution of the above problem exists \textbf{if and only if } $$\mu_1(J\star u -u ,\O)<\lambda<\mu_1(J\star u -u ,\o).$$  Moreover, they have showed that this solution is unique and have established the following asymptotic behaviours:
\begin{align*}
&\lim_{\lambda \to \mu_1(J\star u -u ,\O)}u_\lambda(x) =0 \quad \text{ for all } x \in \O,\\
&\lim_{\lambda \to \mu_1(J\star u -u ,\o)}u_\lambda(x) =+\infty \quad \text{ for all } x \in \O.
\end{align*}

These results have been recently extended to the more general equation \eqref{refuge-eq} with a quadratic nonlinearity ($s(a(x)-b(x)s)$) and under some assumptions on the symmetry of the kernel $K$ and some extra conditions on $a$ and $\lambda$, see \cite{Li2012}.
  
Here we address these questions of existence, uniqueness and the description of some bifurcation diagrams for a general kernel $K$ and with no restriction on the coefficients $a_i,\lambda$  and $\beta$. 
 
 In what follows we will always assume that the functions $a_i$ and $\beta$ satisfy:
 \begin{equation}\label{refuge-hyp1}
 \begin{cases}
&a_i, \beta \in   C(\O), a_1 \ge 0, \beta \ge 0 \\
&\O\setminus \overline{supp(a_1)} \text{ is a open set of } \R^n\\
&K \in  C(\R^n \times \R^n ), K \ge 0
\end{cases}
 \end{equation}
 For the dispersal kernel, we will also require that $K$ satisfies:
 \begin{equation}\label{refuge-hyp2}
  \exists c_0 > 0,  \eps_0 > 0 \quad \text{such that} 
  \inf_{ x\in \O}\left( \inf_ {y\in B(x,\eps_0 )} K(x,y)\right) > c_0.
 \end{equation}
 
 A typical example of such dispersal kernel is given by
$$K(x,y)=J\left(\frac{x_1-y_1}{g_1(y)h_1(x)};\frac{x_2-y_2}{g_2(y)h_2(x)};\ldots; \frac{x_n-y_n}{g_n(y)h_n(x)}\right),$$ 
with $J\in C(\R^n)$ continuous, $J(0)>0$  and  $0<\alpha_i\le g_i\le \beta_i$ and $0\le h_i\le \beta_i$. Such type of kernel have been recently introduced in \cite{Cortazar2007a} to model a nonlocal heterogeneous dispersal process. 
To simplify the presentation of our results, we also introduce the notation $\oplb{u}{\O}$ for the continuous linear operator 
$$\oplb{u}{\O}:=\int_{\O}K(x,y)u(y)\,dy -k(x)u(x).$$

In \cite{Garcia-Melian2009a,Li2012} the analysis essentially relies on the existence of positive eigenfunction associated with a principal eigenvalue $\mu_1$ and a $L^2$  variational characterisation of $\mu_1$. However, such properties ( existence of a positive eigenfunction and a $L^2$ variational characterisation of $\mu_1$)   does not hold for general kernels $K$ and $a_i$ \cite{Coville2013c} and a new approach and characterisation of the principal eigenvalue has to be developed.

In the past few years, the spectral properties of nonlocal operators such as $\lb{\O}+a$  have been intensively studied  \cite{Bates2007,Coville2010,Coville2013c,Coville2008,Coville2012,Garcia-Melian2009,Hutson2003,Hutson2008,Kao2010}. In particular a notion of generalized principal eigenvalue $\mu_p$ of a linear operator  $\lb{\O}+a$ has been introduced in \cite{Coville2010, Coville2012} and is defined by 
$$ \mu_p(\lb{\O}+a):=\sup\{\mu\in \R\,|\,\exists \phi\in C(\bar \O), \phi>0, \text{ so that } \; \oplb{\phi}{\O}+(a+\lambda)\phi \le 0\}.$$
 $\mu_p$ is called a generalized principal eigenvalue because $\mu_p$ is not necessarily associated with a $L^1$ positive eigenfunction \cite{Coville2010,Coville2013c,Kao2010,Shen2012a}. Such notion has been successfully used to derive an optimal criterium for the existence of a unique positive solution of \eqref{refuge-eq} in absence of a refuge zone \cite{Coville2010,Coville2012}.     

Equipped with this notion of generalised eigenvalue, we can now state our results.  
We first present an optimal criterium for the existence of a unique positive bounded solution to \eqref{refuge-eq}. Namely, we show 
\begin{theorem}\label{refuge-thm1}
Let $K$, $a_i$, $\beta$ satisfy the assumptions \eqref{refuge-hyp1}--\eqref{refuge-hyp2} and let $\o$ be the refuge set
$$\o:=\{x\in \bar \O| \beta(x)=0\}.$$
Then a positive continuous bounded solution $u$ of \eqref{refuge-eq} exists \textbf{ if and only if }
$$\mu_p(\lb{\O}+a_0+\lambda a_1)<0<\mu_p(\lb{\o}+a_0+\lambda a_1),  $$
where we set $\mu_p(\lb{\o}+a_0+\lambda a_1)=-\sup_{\o}(a_0+\lambda a_1)$ when  $\stackrel{\circ}{\o}=\emptyset$.
Moreover the solution is unique.
\end{theorem}

 Next we analyse the partially controlled problem  \eqref{refuge-eq} i.e. we describe the set  $\{\lambda,u_\lambda\}$ where $u_\lambda$ is a positive bounded continuous solution to \eqref{refuge-eq}. 
 We start by  describing $\{\lambda,u_\lambda\}$ in a case of the absence of a refuge zone. We prove the following
 \begin{theorem}\label{refuge-thm2}
 Assume that $K,a_i$ and $\beta$ satisfy \eqref{refuge-hyp1}--\eqref{refuge-hyp2}. Assume further that $\beta>0$ in $\bar \O$ then there exists $\lambda^* \in [-\infty;\infty)$, so that for all $\lambda > \lambda^*$ there exists a unique positive continuous solution $u_\lambda$ to \eqref{refuge-eq}. When $\lambda^*\in \R$, there is no positive solution to \eqref{refuge-eq} for all $\lambda \le \lambda^*$. Moreover, we have the following 
 trichotomy: 
 \begin{itemize}
 \item $\lambda^*=-  \infty$ when $\mu_p(\lb{\O\setminus \O_1}+a_0)<0$,
 \item $\lambda^*\in [-  \infty,\infty)$ when $\mu_p(\lb{\O\setminus \O_1}+a_0)=0$.
 \item $\lambda^*\in \R$ when $\mu_p(\lb{\O\setminus \O_1}+a_0)>0$.
 \end{itemize}
 In addition, the map $\lambda\to u_\lambda$ is monotone increasing and we have
\begin{align*} 
& \forall\, x \in \bar \O\quad \lim_{\lambda \to +\infty} u_\lambda (x) = +\infty,\\
&\forall\, x \in \bar \O\quad \lim_{\lambda \to \lambda^{*,+}} u_\lambda(x) = u_\infty(x),
\end{align*}
 where $u_\infty \equiv 0$ on $\O_1 := \{x\in \O| \,a_1 (x) > 0\}$ and $u_\infty$ is a nonnegative solution to 
 $$\int_{\O\setminus \O_1} K(x,y)u(y)dy - k(x)u + a_0 (x)u -\beta u^p = 0 \quad \text{in}\quad \O\setminus\O_1.$$
 Furthermore, $u_\infty$ is non trivial when $\mu_p(\lb{\O\setminus \O_1}+a_0)<0$.
 \end{theorem}

Finally, we describe the set $\{\lambda,u_\lambda\}$ in the situation where a refuge zone exists.  We  prove the following

\begin{theorem}\label{refuge-thm3}
 Assume that $K,a_i$ and $\beta$ satisfy \eqref{refuge-hyp1}--\eqref{refuge-hyp2}. Assume further that $\o\neq\emptyset,$  then there exists two quantities $\lambda^*,\lambda^{**} \in [-\infty,+\infty]$ so that  we have the following dichotomy :
\begin{itemize}
\item Either $\lambda^{**}\le \lambda^{*}$ and there exists no positive bounded  solution to  \eqref{refuge-eq}.
\item Or  $\lambda^*< \lambda^{**}$, and  for all $\lambda \in (\lambda^*,\lambda^{**})$ 
 there exists a unique positive bounded continuous solution  to \eqref{refuge-eq}. When $\lambda^*, \lambda^{**}\in \R$, there is no positive bounded solution to \eqref{refuge-eq} for all $\lambda \le \lambda^*$ and for all $\lambda\ge \lambda^{**}$. Moreover, the map $\lambda\to u_\lambda$ is monotone increasing and we have
\begin{itemize}
\item[(i)] 
$$ \lim_{\lambda \to \lambda^{**,-}} \|u_\lambda\|_{\infty,\o}= +\infty, $$
where $\|u_\lambda\|_{\infty,\o}:=\sup_{x\in\o}|u_\lambda(x)|$.
\item[(ii)] If $\mu_p(\lb{\o}+a_0+\lambda^{**}a_1)$ is an eigenvalue in $L^1(\o)$ or $\lambda^{**}=+\infty$ then  $$ \forall\, x \in \bar \O, \lim_{\lambda \to \lambda^{**,-}} u_\lambda (x) = +\infty.$$
\item[(iii)] For all $x \in \bar \O$ we have $ \lim_{\lambda \to \lambda^{*,+}} u_\lambda(x) = u_\infty(x),$
where $u_\infty$ is a  function satisfying on $\O_1 := \{x\in \O| \,a_1 (x) > 0\, u_\infty\equiv 0$ and $u_\infty$ is a nonnegative solution to 
 $$\int_{\O\setminus \O_1} K(x,y)u(y)dy - k(x)u + a_0 (x)u -\beta u^p = 0 \quad \text{in}\quad \O\setminus\O_1.$$
 Furthermore, $u_\infty$ is non trivial when $\mu_p(\lb{\O\setminus \O_1}+a_0)<0$.
\end{itemize}
 \end{itemize}
\end{theorem}

Before going to the proofs of theses results we would like to make some additional comments.
The assumption $ $ can be relaxed and we can get a full description of the curves when $a_1>0$ in $\bar \O$.

The paper is organised as follows. In a preliminary  section we recall some known results on $\mu_p$ and on the positive solution of a KPP equation. Then in Section \ref{refuge-section-iff} we prove the existence criterium of Theorem \ref{refuge-thm1}. The proof of Theorem \ref{refuge-thm2} is done in  Section \ref{refuge-section-kpp}.  Finally, in Section \ref{refuge-section-refuge} we analyse the bifurcation diagram of \eqref{refuge-eq} in the presence of a refuge zone (Theorem \ref{refuge-thm3}).


 \section{Preliminaries}
 In this section, we recall some results on the principal eigenvalue of a linear nonlocal operator and some known results about the KPP equation below
 \begin{equation}\label{refuge-eq-kpp}
 \oplb{u}{\O}+f(x,u)=0 \quad \text{ in }\quad \O
 \end{equation}
 where $$\oplb{u}{\O}:=\int_{\O}K(x,y)u(y)\,dy -k(x)u(x)$$
 and $f(x,s)$ is satisfying
\begin{align}
\label{refuge-hyp-kpp} \left\{
\begin{aligned}
& \hbox{$f \in C(\O\times[0,\infty))$ and is differentiable with respect to $s$} \\
& \hbox{$f_u(\cdot,0) \in C(\O)$}  \\
& \hbox{$f(\cdot,0)\equiv 0$ and $f(x,s)/s$ is decreasing with
respect to $s$}
\\
& \hbox{there exists $M>0$ such that $f(x,s)\le 0$ for all $s \ge
M$ and all $x$.}
\end{aligned}
\right.
\end{align}
The simplest example of such a nonlinearity is
$$
f(x,u) = u (a(x) - u ),
$$
where $a(x)\in C(\O)$.

It has been shown in \cite{Bates2007,Coville2010} that the existence of a positive solution of \eqref{refuge-eq-kpp} is conditioned to the sign of the principal eigenvalue $\mu_p$ of the linear operator $\lb{\O}+f_u (x,0)$  where $\mu_p$  is defined by the formula
 $$\mu_p (\lb{\O} + f_u (x,0)) := \sup \{\mu \in \R \, |\, \exists \, \phi \in C(\O), \phi > 0 \;\text{ so that }\; \oplb{\phi}{\O}+ f_u (x,0)\phi + \mu\phi\le 0\}.$$
 That is to say
 \begin{theorem}[\cite{Bates2007,Coville2010}]\label{refuge-thm-existence} Let $\O$ be a bounded open set anf assume that $K$ and $f$ satisfy respectively \eqref{refuge-hyp1}--\eqref{refuge-hyp2} and \eqref{refuge-hyp-kpp}. Then there exists a unique positive continuous solution  to \eqref{refuge-eq-kpp} \textbf{ if and only if} $\mu_p(\lb{\O}+ f_u (x,0)) < 0$. Moreover,  if   $\mu_p \ge 0$ then any  non negative uniformly bounded solution  of \eqref{refuge-eq-kpp} is identically zero.
 \end{theorem}
 Also noted in \cite{Coville2010} the principal eigenvalue is not always achieved. This means that there is not always a positive continuous eigenfunction associated with $\mu_p$. However as shown in \cite{Coville2013c}, we can  always  associate a  positive measure $d \mu$  with $\mu_p$. 
 More precisely,
 
 \begin{theorem}[\cite{Coville2013c}]\label{refuge-thm-mu}
Let $\O$ be an open bounded set and assume that  $K$ and $f_u(x,0)$ satisfy  the assumptions \eqref{refuge-hyp1} and \eqref{refuge-hyp-kpp}. Then there exists a positive measure  $d\mu \in \M^+(\O)$, so that for any $\phi\in C_c(\O)$ we have 
$$\int_{\O}\phi(x)\left(\int_{\O}K(x,y)d\mu(y)\right) dx +\int_{\O}\phi(x)( f_u(x,0)-k(x)+\mu_p)d\mu(x)=0.$$
Moreover,  there exists a positive function $\phi_p \in L^1(\O)\cap C(\O\setminus \Sigma)$ so that $\inf_{\O}\phi_p>0$ and  $d\mu(x)= \phi_p(x) dx + d\mu_s(x)$ where $d\mu_s(x)$ is a non negative singular measure with respect to the Lebesgue measure whose support lies in the set $\Sigma:=\{y\in \bar \O|f_u(y,0)-k(y)=\sup_{x\in \O}(f_u(x,0)-k(x))\}$.  
\end{theorem}

 As proved in \cite{Coville2010,Kao2010,Shen2012a}, when $\O$ is an open bounded set we can find a condition on the coefficients which guarantees that $d\mu_s(x)\equiv 0$ and the existence of a positive continuous eigenfunction. For example the existence of principal eigenfunction is guaranteed,  if we assume that the function $a(x) := f_u (x,0) -\int_{\O}K(y,x)dy$ satisfies
 $$\frac{1}{\sup_{\O}a -a(x)}\not\in L^1 (\O' ) \quad\text{ for some open bounded domain }\quad \O'\subset \bar\O.  $$
 
For the existence of principal eigenfunction as remark in \cite{Coville2012} we also have this useful criteria:
\begin{proposition} \label{refuge-prop-phip}
 Let $\O$ be a bounded open set, then there exists a positive continuous eigenfunction associated to 
$\mu_p$ \textbf{if and only if} $\mu_p(\lb{\O}+a)< -\sup_{\O}a$. 
 \end{proposition}
 
 Next we recall some properties of the principal eigenvalue $\mu_p$ that we will constantly use along this paper:
 
 \begin{proposition}\label{refuge-prop-mup}
\begin{itemize}
\item[(i)] Assume $\O_1\subset\O_2$, then for the two operators
\begin{align*}
&\oplb{u}{\O_1}+a(x)u:=\int_{\O_1}K(x,y)u(y)\,dy-k(x)u +a(x)u\\
&\oplb{u}{\O_2}+a(x)u:=\int_{\O_2}K(x,y)u(y)\,dy-k(x)u +a(x)u
\end{align*}
respectively defined on $C(\O_1)$ and $C(\O_2)$ we have
 $$
\mu_p(\lb{\O_1}+a(x))\ge \mu_p(\lb{\O_2}+a(x)).
$$
\item[(ii)]Fix $\O$  and assume that $a_1(x)\ge a_2(x)$, then 
$$
\mu_p(\lb{\O}+a_2(x))\ge\mu_p(\lb{\O}+a_1(x)).
$$
Moreover, if $a_1(x)\ge a_2(x)+\delta$ for some $\delta>0$ then 
$$
\mu_p(\lb{\O}+a_2(x))> \mu_p(\lb{\O}+a_1(x)).
$$
\item[(iii)] $\mu_p(\lb{\O}+a(x))$ is Lipschitz continuous in $a(x)$. More precisely,
$$|\mu_p(\lb{\O}+a(x))- \mu_p(\lb{\O}+b(x))|\le \|a(x)-b(x)\|_{\infty}$$

\item[(iv)] Assume $\O_1\subset\O_2$, then for the two operators
\begin{align*}
&\oplb{u}{\O_1}+a(x)u:=\int_{\O_1}K(x,y)u(y)\,dy-k(x)u +a(x)u\\
&\oplb{u}{\O_2}+a(x)u:=\int_{\O_2}K(x,y)u(y)\,dy-k(x)u +a(x)u
\end{align*}
respectively defined on $C(\O_1)$ and $C(\O_2)$. Assume  that the corresponding principal eigenvalue are associated to a positive continuous principal eigenfunction. Then we have 
$$|\mu_p(\lb{\O_1}+a(x))- \mu_p(\lb{\O_2}+a(x))|\le C_0|\O_2\setminus \O_1|,$$
where $C_0$ depends on $K$ and $\phi_2$.
\item[(v)]  We always have the following estimate
$$-\sup_{\O}\left(a(x)+\int_{\O}K(x,y)\,dy\right)\le \mu_p(\lb{\O}+a)\le -\sup_{\O}a.$$
\end{itemize}
\end{proposition}
\dem{Proof:}
We refer to \cite{Coville2010} for the proofs of $(i)-(iii)$  and $(v)$, so we will be concerned only with $(iv)$. Let us introduce the following quantity:
$$\mu_p' (\lb{\O} + a) := \inf \{\mu \in \R \, |\, \exists \, \phi \in C(\bar \O), \phi > 0 \;\text{ so that }\; \oplb{\phi}{\O}+ a\phi + \mu\phi\ge 0\}.$$
 One can check that $\mu_p (\lb{\O} + a) = \mu'_p (\lb{\O} + a)$. Indeed since there is a positive eigenfunction associated with $\mu_p (\lb{\O} +a)$  one has $\mu'_p (\lb{\O} +a) \le \mu_p (\lb{\O} +a)$ by definition of $\mu'_p (\lb{\O} +a)$. We obtain the equality by arguing as follows. Assume by contradiction that $\mu_p (\lb{\O} + a) > \mu'_p (\lb{\O} + a)$. Then there exists $\mu$ so that $\mu'_p (\lb{\O} + a) < \mu < \mu_p (\lb{\O} + a)$ and from the definition of $\mu_p$ and $\mu'_p$ there exists two  positive continuous functions $\psi$ and $\phi$ so that
\begin{align*}  
  &\oplb{\phi}{\O}  + a(x)\phi + \mu\phi \ge 0,\\
  &\oplb{\psi}{\O}  + a(x)\psi + \mu\psi < 0.
  \end{align*}
  From the last inequalities we deduce that $\psi > 0$ in $\bar \O$ and by setting $w :=\frac{\phi}{\psi}$ it follows that
\begin{align*} 
 0 \le \oplb{\phi}{\O}  + (a(x) + \mu)\phi &= \oplb{\phi}{\O}  + (a(x) + \phi)\frac{\phi}{\psi}\psi,\\
 &\le \oplb{\phi}{\O}  - \frac{\phi}{\psi} (x)\oplb{\psi}{\O},\\
 &\le \int_{\O} K(x,y)\psi(y)(w(y) - w(x))\,dy.
 \end{align*}
 Thus $w$ cannot achieve a maximum in $\bar \O$ without being constant. $w$ being continuous in $\bar \O,$  it follows  that $\phi = c\psi$ for some positive constant $c$. Thus we get the contradiction
 $$0 \le \oplb{\phi}{\O}  + (a(x) + \mu)\phi = c\left( \oplb{\psi}{\O}  + (a(x) + \mu)\psi\right) < 0.$$ 
We are now in position to prove (iv). Let $\phi_2$ be the eigenfunction associated to$ \mu_p (\lb{\O_2} + a(x))$ normalized by 
$\|\phi_2\|_{\infty} = 1$ and let us set $C_0 := \frac{\| K(\cdot,\cdot)\|_{\infty}}{\min_{\O_2}\phi_2}.$  

Now, let us  show that $(\phi_2 ,\mu_p (\lb{\O_2} + a) + C_0 |\O_2\setminus\O_1 |)$ is an adequate test function for $\mu'_p (\lb{\O_1} + a)$. By a direct computation and by using the normalisation of $\phi_2$ we have
\begin{align*}
 \oplb{\phi_2}{\O_1} + (a + \mu_p (\lb{\O_2} + a) + C_0 |\O_2\setminus\O_1 |))\phi_2 &= -\int_{\O_2\setminus\O_1}K(x,y)\phi_2 (y)\,dy + C_0 |\O_2\setminus\O_1 |\phi_2,\\
 &\ge -\| K(\cdot,\cdot)\|_{\infty}|\O_2\setminus\O_1 | + C_0 |\O_2\setminus\O_1 |\phi_2,\\
 &\ge  \left( \frac{\phi_2}{\min_{\O_2} \phi_2} - 1\right)\| K(\cdot,\cdot)\|_{\infty}|\O_2\setminus\O_1 |.
 \end{align*}
 Therefore $\mu'_p (\lb{\O_1} + a) \le \mu_p (\lb{\O_2} + a) + C_0 |\O_2\setminus\O_1 |,$
 which combined with $(i)$ and using that $\mu'_p (\lb{\O_1} + a) = \mu_p (\lb{\O_1} + a)$ leads to
 $$|\mu_p (\lb{\O_1} + a(x)) - \mu_p (\lb{\O_2} + a(x))| \le C_0 |\O_2\setminus\O_1 |.$$
\fdem

\section{Optimal existence criterium} \label{refuge-section-iff}

 In this section we establish an optimal criterium for the existence of a positive continuous bounded solution to 
 \begin{equation}\label{refuge-eq-refuge}
  \oplb{u}{\O}  + a(x)u -\beta(x)u^p = 0\quad \text{ in }\quad \O,
  \end{equation}
 when there exists $\o\subset\O$ so that $\beta_{|\o}\equiv 0$. Note that \eqref{refuge-eq-refuge} is a particular case of \eqref{refuge-eq-kpp} with $f(x,s) := a(x)s - \beta(x)s^p.$  However, due to the presence of refuge zone (i.e. $\beta_{|\o}\equiv 0$) the function $f(x,s)$ does not satisfy the assumptions \eqref{refuge-hyp-kpp} and the Theorem \ref{refuge-thm-existence} does not apply. But we still have a complete characterisation of the existence of a bounded positive solution. Namely we can show the following Theorem:
 \begin{theorem}\label{refuge-thm-refuge}
 Assume that $K,a$ and $\beta$ satisfy \eqref{refuge-hyp1}--\eqref{refuge-hyp2}. Assume further that there exists $\o\subset \O$ so that $\beta_{|\o}\equiv 0$. Then there exists a bounded positive continuous solution to \eqref{refuge-eq-refuge} \textbf{if and only if} 
 $$\mu_p (\lb{\o} + a) > 0 > \mu_p (\lb{\O} + a),$$
 where we set $\mu_p(\lb{\o}+a)=-\sup_{\o}a$ when  $\stackrel{\circ}{\o}=\emptyset$.
 \end{theorem}
\dem{Proof:} 

First let us assume that $\mu_p (\lb{\o} + a) \le 0$, we will show that there is no positive bounded solution to \eqref{refuge-eq-refuge}. Let us suppose by contradiction that there exists $u$, a positive bounded solution to \eqref{refuge-eq-refuge}. So in $\o, u$ satisfies $$\oplb{u}{\O} + au = 0,$$ which implies that $\max_{\bar \o}  a < 0$ and $u$ is continuous on $\bar \o.$ Furthermore, we have
\begin{equation}\label{refuge-eq-noex1}
  \oplb{u}{\o}  + au = -\int_{\O\setminus\o}K(x,y)u(y)\,dy \le 0.
 \end{equation}
 If $\stackrel{\circ}{\o}=\emptyset$ then we obtain easily a contradiction. Indeed in such case, we have $\mu_p(\lb{\o}+a)=-\sup_{\o}a$ which leads to the contradiction
 $$0<-\sup_{\o}a=\mu_p(\lb{\o}+a)\le 0.$$
 In the other situations,  $\stackrel{\circ}{\o}\neq\emptyset$ and to obtain our desired contradiction we argue as follows. Since $\mu_p (\lb{\o} + a) \le 0 < -\max_{\bar \o} a$, by  Proposition \ref{refuge-prop-phip} there exists a positive continuous eigenfunction associated with $\mu_p (\lb{\o} + a).$ As a consequence there exists also a positive continuous eigenfunction associated with $\mu_p (\lbs{\o} + a)$ where $ \lbs{\o}+ a$  is defined by
 $$
 \oplbs{\phi}{\o} + a\phi := \int_{\o} K(y,x)\phi(y)dy -k(x)\phi + a(x)\phi.$$
 We can easily check that $\mu_p (\lb{\o} + a) = \mu_p (\lbs{\o} + a).$ Let us denote by $\phi^*$ the positive continuous principal eigenfunction associated with $\mu_p (\lbs{\o} + a).$ Now by multiplying \eqref{refuge-eq-noex1} by $\phi^*$ and then integrating  over $\o$,  it  follows 
$$
\int_{\o} \phi^*(x)\oplb{u}{\o}(x)dx + au\phi^* (x)dx\le - c_0\int_{\o}\phi^*\left(\int_{\O\setminus\o}K(x,y)\, dy \right).$$ 

 By using Fubini's Theorem in the above inequality we get  the contradiction
 $$0 \le -\mu_p (\lbs{\o} + a)\int_{\o}\phi^* u\le -c_0 \int_{\o}\phi^*\left(\int_{\O\setminus\o}K(x,y)\, dy \right)<0.$$
Thus in both cases, there is no bounded solution to \eqref{refuge-eq-refuge} when   $\mu_p (\lb{\o} + a) \le 0.$
 
 Next we see that there is no positive bounded solution for \eqref{refuge-eq-refuge} when $\mu_p (\lb{\O} + a) \ge 0.$ In this situation, with some modifications we can reproduce the argumentation developed in \cite{Coville2010} (Subsection 6.2).  Let us assume that a positive solution of \eqref{refuge-eq-refuge} exists and let us denote $u$ this solution. 
 We first observe that by following the argument developed in \cite{Bates2007} we can see that $u$ is continuous in $\O$ and  there exists positive constants $\delta$ and $c_0$ so that 
 $$\begin{cases}
&\inf_{\O} u \ge c_0,\\
 &\inf_{x\in \O} (k(x) - a(x) + \beta(x)u^{p-1} ) \ge \delta.
 \end{cases}
 $$
 From the monotone behaviour with respect to the $s$  of the function $g(x,s) := (a -\beta(x)s^{p-1} ),$ we deduce that $a -\beta u^{p-1}\le  a(x) - \beta c_0^{p-1} \le a.$ Now let us denote $\gamma(x) = a(x) - \beta(x) c_0^{p-1}.$ By construction, we have $\gamma(x) \le a(x)$ and we see by $(ii)$ of Proposition \ref{refuge-prop-mup} that
 $$\mu_p (\lb{\O} + \gamma(x)) \ge \mu_p (\lb{\O} + a(x)) \ge 0.$$
 Moreover, since $u$ is a solution of \eqref{refuge-eq-refuge}, we have
 \begin{equation}\label{refuge-eq-noex2}
 \oplb{u}{\O}  + \gamma u \ge \oplb{u}{\O} + au - \beta u^p = 0,
 \end{equation}
with a strict inequality for any $x \in\O\setminus\o.$ 
 
 We claim that
\begin{cla} 
 There exists $\delta > 0$ and a positive continuous function $\phi$ so that $\inf_{\O}\phi> \delta$ and
 $$\oplb{\phi}{\O}  + \gamma\phi \le 0.$$
 \end{cla}
 Assume for the moment that the Claim holds true then we get our desired contradiction by arguing as follow. Since $\phi > \delta$ we can define the following quantity
 $$\tau^* := \inf\{ \tau > 0 | u \le \tau \phi\}.$$
 Obviously, by proving that $\tau^* = 0$ we get the contradiction
 $$c_0 \le u \le 0.$$
 Assume by contradiction that $ \tau^* > 0$ and let us denote $w:=\tau^*\phi -u$. By definition of $\tau^*$, there exists $x_0 \in\bar \O$ such that $\tau^*\phi(x_0 ) = u(x_0 ) > 0$ and from \eqref{refuge-eq-noex2} we see that  $w$ satisfies
  $$\oplb{w}{\O}  + \gamma w \le 0.$$
  By evaluating the above expression at $x_0$, since $w\ge 0$ we see that  
$$0\le \int_{\O}K(x_0,y)w(y)\,dy  \le 0.$$
 Therefore, since K satisfies \eqref{refuge-hyp2} we must have $w(y) = 0$ for almost every $y \in \bar \O.$ Thus, we end up with $\tau^*\phi\equiv u$ and we get the following contradiction
 $$ 0 < \oplb{u}{\O}  + \gamma u = \oplb{\tau^*\phi}{\O}  + \gamma \tau^*\phi \le  0 \quad \text{ on }\quad \O\setminus\o.$$
 Hence $\tau^* = 0.$ 
 \dem{Proof of the Claim:} 
When $\mu_p (\lb{\O}+\gamma) > 0$ then by definition of the principal eigenvalue for all positive $0 < \mu <\mu_p (\lb{\O}+\gamma)$ there exists a positive continuous function $\phi$ such that
 $$\oplb{\phi}{\O}+\gamma\phi \le -\mu\phi < 0.$$
 Observe that $\phi\ge \delta$ for some positive $\delta$ since otherwise there exists $x_0\in \bar \O$  so that $\phi(x_0 ) = 0$ and we get the contradiction
 $$0<  \oplb{\phi}{\O}(x_0 ) + \gamma(x_0 ) \phi(x_0) \le 0.$$
 When $\mu_p (\lb{\O}+\gamma) = 0$ we argue as follows. By construction,  $a \ge \gamma$ and 
  on $\bar \O\setminus \o$  we have 
  \begin{equation}\label{refuge-eq-noex3}
  \gamma < a \le \sup_{\O} a\le -\mu_p(\lb{\O}+a)\le 0.
  \end{equation} 
  
 And another hand on $\o$ since $\beta_{|\o} \equiv 0$, we  have
 $$\oplb{u}{\O} + a(x)u = 0,$$
 which leads to $\sup_{\bar \o}  a < 0$. So on $\bar \o$ we also have
 \begin{equation}
 \gamma \le  \sup_{\o} a < 0. \label{refuge-eq-noex4}
 \end{equation}
 By combining \eqref{refuge-eq-noex3} and  \eqref{refuge-eq-noex4} it follows that $\sup_{\bar\O} \gamma < 0.$ 
 Now, since  $ 0 = \mu_p (\lb{\O}+\gamma) < -\sup_{\bar \O} \gamma$ we deduce from  Proposition \ref{refuge-prop-phip} that there exists a continuous positive principal eigenfunction $\phi$ associated with $\mu_p (\lb{\O}+\gamma).$ 
 As above, we have $\inf_{\O} \phi > \delta$ for some positive $\delta$. 
\fdem
 
 Lastly, let us construct a positive bounded solution to \eqref{refuge-eq-refuge} when the condition
\begin{equation}
 \mu_p (\lb{\o} + a) > 0 > \mu_p (\lb{\O} + a)\label{refuge-eq-cond}
 \end{equation}
 is satisfied. The uniqueness of this solution follows form a similar argumentation as in \cite{Bates2007,Coville2010}, so we will omit the proof here. 
 
 From the condition $0 > \mu_p (\lb{\O} + a)$, by reproducing the argument in \cite{Coville2010} we can find a positive bounded subsolution $\phi_0$ of the problem \eqref{refuge-eq-refuge} so that $\kappa\phi_0$ is  still a subsolution for any $\kappa$ small and positive. Here the main difficulty is to find a positive supersolution $\psi$. Indeed,  due to the existence of  a refuge zone, the large positive constants are not supersolutions of \eqref{refuge-eq-refuge}.
We claim 
\begin{cla}\label{refuge-cla-super}
When the condition \eqref{refuge-eq-cond} is satisfied, then there exists $\psi>0$, $\psi\in C(\bar \O)$ supersolution of \eqref{refuge-eq-refuge}
\end{cla}

Note that by proving the claim we end the construction of the solution to \eqref{refuge-eq-refuge}. Indeed, since for $\kappa$ small we have $\kappa \phi \le \psi$,  by the monotone iterative scheme there exists a solution $u$ to \eqref{refuge-eq-refuge} so that $\kappa \phi \le u\le \psi$.

 \fdem
 
Now, let us turn our attention to the proof of the Claim.
\dem{Proof of the Claim:}
Let us first assume that $\stackrel{\circ}{\o}\neq \emptyset$.

In this situation, by following the argument in \cite{Coville2010} (Subsection 6.1) we can introduce a regularisation $a_\eps \in C(\O)$ of $a -k$ so that the following operator
 $$\oplb{u}{\eps,\o}:= \int_{\o}K(x,y)u(y)dy + a_\eps (x)u$$
 has a positive continuous principal eigenfunction. 
 By  continuity of $\mu_p (\lb{\eps,\o})$ with respect to $a_\eps$ ($(iii)$ of Proposition \ref{refuge-prop-mup}) we can find $\eps$ small so that
 $$\mu_p (\lb{\eps,\o} ) -\|a_\eps - a + k\|_{\infty}\ge \frac{\mu_p (\lb{\o} + a)}{2}.$$
 
 Let $\eps$ be fixed and let us denote $\o_\delta$ the following set
$$\o_\delta := \{x \in \O\;|\;  d(x;\o) < \delta\}.$$
 
By  continuity of the function $\sup_{\o_\delta} a_\eps$ with respect to $\delta$, there exists $\delta_0$ so that for all $\delta \le \delta_0$ we have 
$$|\sup_{\o_\delta}a_\eps - \sup_{\o}a_\eps|\le\frac{-\mu_p(\lb{\o}+a_\eps) -\sup_{\o}a_\eps}{2}.  $$
 So, by $(i)$ of Proposition \ref{refuge-prop-mup}, we deduce from the above inequality that  we have  for all $\delta\le \delta_0$,
 $$\mu_p(\lb{\o_\delta}+a_\eps)\le \mu_p(\lb{\o}+a_\eps)< -\sup_{\o_\delta}a_\eps.$$
 Therefore, thanks to Proposition \ref{refuge-prop-phip} for all $\delta \le \delta_0$  there exists a positive continuous eigenfunction associated with  $\mu_p(\lb{\eps,\o_\delta})$.
 
 By continuity of  $\mu_p(\lb{\eps,\o_\delta})$ with respect to the domain ((iv) of Proposition \ref{refuge-prop-mup}) we achieve for $\delta$ small enough, say $\delta \le  \delta_1,$ 
 \begin{equation}\label{refuge-eq-exis1}
 \mu_p (\lb{\eps,\o} ) \ge \mu_p (\lb{\eps,\o_\delta} ) \ge \mu_p (\lb{\eps,\o_{\delta_1}} ) \ge \|a_\eps - a+k\|_{\infty} +\frac{ \mu_p (\lb{\o} + a)}{ 8}.
 \end{equation}

By construction $\bar \O\setminus \o_ \delta$ and $\bar \o_{\frac{\delta} {2}}$ are two disjoints bounded closed set, so by the Urysohn Lemma there exists a nonnegative continuous function $\eta_1$ such that $0\le \eta_1\le 1, \eta_1 (x) = 1$ in $\bar \O\setminus \o_ \delta , \eta_1 (x) = 0$ in $\bar\o_{\frac{\delta}{ 2}}.$

Let $\psi_1,\psi_2$ be the  following continuous functions
 $$\psi_1 :=\begin{cases} 
 &C_1 \eta_1 \quad\text{ in }\quad \O\setminus\o_{\frac{\delta}{ 2}}\\
 &0 \quad\text{elsewhere},
\end{cases} 
\qquad
 \psi_2 :=\begin{cases} 
&C_2(1-\eta_1)\Psi_\delta \quad \text{ in }\quad \o_\delta\\
 &0 \quad\text{elsewhere}.
 \end{cases}
 $$
 where $\Psi_\delta$ denotes the positive continuous eigenfunction associated with  $\mu_p (\lb{\eps,\o_\delta} )$ normalized by $\|\Psi_\delta\|_{\infty}= 1$ and $C_1$ and $C_2$ are positive constants to be specified later on. 
 Consider now the function $\psi:= \sup (\psi_1 ,\psi_2 )$, we will prove that for well chosen $C_1$ and $C_2$, $\psi$ is a supersolution of \eqref{refuge-eq-refuge}. 
 
 On $\O\setminus \o_{\frac{\delta}{2}}$,  a short computation shows that for $C_1$ large 
\begin{align*}
 \oplb{\psi}{\O}  + a\psi -\beta \psi^p &\le C_1^{\frac{p+1}{2}}\left(\int_{\O} K(x,y)\,dy + a - \beta C_1^{\frac{p-1}{2}}\right),\\ 
&\le C_1^{\frac{p+1}{2}}\left(\int_{\O} K(x,y)\,dy + a - \inf_{\O\setminus \o_{\frac{\delta}{2}}}(\beta) C_1^{\frac{p-1}{2}}\right). 
\end{align*}

By construction $ \inf_{\O\setminus \o_{\frac{\delta}{2}}}(\beta)>0$ and  $\frac{p-1}{2}>0$, so for $C_1$ large enough we have on $\O\setminus \o_{\frac{\delta}{2}}$ 
\begin{equation}\label{refuge-eq-exis3}
 \oplb{\psi}{\O}  + a\psi -\beta \psi^p  \le 0.
\end{equation}

Now, observe that for $C_2\ge   C_1$ large, on $\o_{\frac{\delta}{2}}$ we have
\begin{multline*}
 \oplb{\psi}{\O}  + a\psi -\beta \psi^p \le C_1\int_{\O\setminus\o_{\delta}} K(x,y)\,dy+C_2\int_{\o_\delta\setminus \o_{\frac{\delta}{2}}}K(x,y)\,dy \\+C_2\left(\int_{\o_{\frac{\delta}{2}}}K(x,y)\Psi_\delta(y)\,dy -k(x)\Psi_\delta +a(x)\Psi_\delta\right).
 \end{multline*}

Since $\Psi_\delta>0$ in $\o_\delta$, we have 
\begin{multline*}
 \oplb{\psi}{\O}  + a\psi -\beta \psi^p \le C_2\left (\frac{ C_1}{C_2}\bar K_1+\bar K_2|\o_\delta \setminus \o_{\frac{\delta}{2}}|\right)\\+C_2\left(\int_{\o_\delta}K(x,y)\Psi_\delta(y)\,dy -k(x)\Psi_\delta +a(x))\Psi_\delta\right),
 \end{multline*}
 where $ \bar K_1:=\sup_{x\in\O}\int_{\O}K(x,y)\,dy$ and $ \bar K_2:=\|K(\cdot,\cdot)\|_{\infty}$.
 
 Recall that $\Psi_\delta$ is the eigenfunction associated with $\mu_p(\lb{\eps,\o})$, so it follows  that
$$
 \oplb{\psi}{\O}  + a\psi -\beta \psi^p  \le C_2\left (\frac{ C_1}{C_2}\bar K_1+\bar K_2|\o_\delta \setminus \o_{\frac{\delta}{2}}|+(a-k-a_\eps -\mu_p(\lb{\eps,\o_\delta}))\Psi_\delta\right)
 $$
 
which combined with \eqref{refuge-eq-exis1} reduces to 
$$  \oplb{\psi}{\O}  + a\psi -\beta \psi^p  \le C_2\left (\frac{ C_1}{C_2}\bar K_1+\bar K_2|\o_\delta \setminus \o_{\frac{\delta}{2}}| -\frac{\mu_p(\lb{\o}+a)}{8}\Psi_\delta\right).$$

By using that for all $\delta \in [0,\delta_0 ]$, the principal eigenfunction $\Psi_\delta$ associated to $\mu_p (\lb{\eps,\o_\delta} )$ is positive and continuous, we can see that 
 $$\inf_{\delta \in [0,\delta_0]} \inf_{\o_\delta} \Psi_\delta > c,$$
for some positive constant $c$.  
 Moreover we can find $\delta$ small, say $\delta \le \delta_1$ so that  for all $\delta \le \delta_1$
$$\bar K_2|\o_\delta \setminus \o_{\frac{\delta}{2}}|-\frac{\mu_p(\lb{\o})}{8}\Psi_\delta\le -\frac{\mu_p(\lb{\o}+a)}{16}c. $$ 
Thus for $\delta\le \delta_1$, we achieve on $\o_{\frac{\delta}{2}}$ 
 
 $$  \oplb{\psi}{\O}  + a\psi -\beta \psi^p \le C_2\left (\frac{ C_1}{C_2}\bar K_1-\frac{\mu_p(\lb{\o}+a)}{16}c\right).$$
 Now by choosing $C_2:=C_1^{\frac{p+1}{2}}$  we have $\lim_{C_1\to +\infty}\frac{ C_1}{C_2}=0$ since $p>1$.
 So for $C_1$ large enough,  say $C_1\ge C_1^*:= \left(\frac{32 \bar K_1}{\mu_p(\lb{\o}+a)c}\right)^{\frac{2}{p-1}}$ we achieve on $\o_\delta$
 
 \begin{equation}\label{refuge-eq-exis2}
 \oplb{\psi}{\O}  + a\psi -\beta \psi^p \le -cC_1^{\frac{p+1}{2}}\frac{\mu_p(\lb{\o}+a)}{32} < 0.
\end{equation}
 
Hence from \eqref{refuge-eq-exis3} and \eqref{refuge-eq-exis2} we see that  the function $\psi$ is a positive continuous supersolution of \eqref{refuge-eq-refuge}.

Let us now assume that $\stackrel{\circ}{\o}=\emptyset$. In this situation, we have $0<\mu_p(\lb{\o}+a)=-\sup_{\o}a$.
By continuity of $a$ and $K$ there exists $\delta$ small so that 
$$\sup_{x\in \o_\delta}\int_{\o_\delta}K(x,y)dy\le \frac{\mu_p(\lb{\o}+a)}{2}<-\sup_{\o_\delta}a,$$ 
where as above $\o_\delta:=\{x\in \O|d(x,\o)<\delta\}$.

From the above inequality it follows from  $(v)$ of Proposition \ref{refuge-prop-mup}  that $0< \mu_p(\lb{\o_\delta}+a)$.
Let us consider $\beta_\delta:=\beta \eta_1$, where $\eta_1$ is constructed above, then we have   $\beta_\delta\le \beta$ and $\o_{\frac{\delta}{2}}=\{x\in \O| \beta_\delta(x)=0 \}$. 
By construction  $\stackrel{\circ}{\o_{\frac{\delta}{2}}}\neq\emptyset$ and $0<\mu_p(\lb{\o_\delta}+a)\le \mu_p(\lb{\o_{\frac{\delta}{2}}}+a)$, therefore by using  the above arguments there exists a positive continuous supersolution $\psi$ to 
$$\oplb{u}{\O}+a(x)u-\beta_\delta u^p=0 \quad \text{in}\quad \O.$$
Thanks to $\beta_\delta\le \beta$, we have 
$$\oplb{\psi}{\O}+a(x)\psi-\beta \psi^p\le \oplb{\psi}{\O}+a(x)\psi-\beta_\delta \psi^p\le0 \quad \text{in}\quad \O$$
and $\psi$ is our desired supersolution.
\fdem 


 \section{The partially controlled problem: The KPP case} \label{refuge-section-kpp}
 In this section we analyse the dependence in $\lambda$ of the positive continuous solutions to \eqref{refuge-eq} in absence of  a refuge zone and we prove the Theorem \ref{refuge-thm2} that we recall below.  
 More precisely, we look for positive continuous solution of the partially controlled problem:
 \begin{equation}\label{refuge-eq-steady}
 \int_{\O} K(x,y)u(y)\, dy - k(x)u(x) + a_0 (x)u + \lambda a_1 (x)u -\beta u^p = 0
 \end{equation}
 when $\beta>0$ and $\lambda \in \R$. 

In absence of a refuge zone, we can show that there exists a critical value $\lambda^*$ characterising completely the existence/non existence of a positive stationary solution. More precisely we have,
\begin{theorem}
 Assume that $K,a_i$ and $\beta$ satisfy \eqref{refuge-hyp1}--\eqref{refuge-hyp2}. Assume further that $\beta>0$ in $\bar \O$ then there exists $\lambda^* \in [-\infty;\infty)$, so that for all $\lambda > \lambda^*$ there exists a unique positive continuous solution $u_\lambda$ to \eqref{refuge-eq-steady}. When $\lambda^*\in \R$, there is no positive solution to \eqref{refuge-eq} for all $\lambda \le \lambda^*$. Moreover, we have the following 
 trichotomy: 
 \begin{itemize}
 \item $\lambda^*=-  \infty$ when $\mu_p(\lb{\O\setminus \O_1}+a_0)<0$,
 \item $\lambda^*\in [-  \infty,\infty)$ when $\mu_p(\lb{\O\setminus \O_1}+a_0)=0$.
 \item $\lambda^*\in \R$ when $\mu_p(\lb{\O\setminus \O_1}+a_0)>0$.
 \end{itemize}
 In addition, the map $\lambda\to u_\lambda$ is monotone increasing and we have
\begin{align*} 
& \forall\, x \in \bar \O\quad \lim_{\lambda \to +\infty} u_\lambda (x) = +\infty,\\
&\forall\, x \in \bar \O\quad \lim_{\lambda \to \lambda^{*,+}} u_\lambda(x) = u_\infty(x),
\end{align*}
 where $u_\infty \equiv 0$ on $\O_1 := \{x\in \O| \,a_1 (x) > 0\}$ and $u_\infty$ is a nonnegative solution to 
 $$\int_{\O\setminus \O_1} K(x,y)u(y)dy - k(x)u + a_0 (x)u -\beta u^p = 0 \quad \text{in}\quad \O\setminus\O_1.$$
 Furthermore, $u_\infty$ is non trivial when $\mu_p(\lb{\O\setminus \O_1}+a_0)<0$.
 \end{theorem}

 \dem{Proof:}
 In absence of a refuge zone, we observe that the  problem \eqref{refuge-eq-steady} is a particular case of the KPP equation \eqref{refuge-eq-kpp} where the nonlinearity $f$ is given by $f(x,s) := a_0 s+\lambda a_1 s-\beta s^p$. Therefore by the Theorem \ref{refuge-thm-existence}, for each $\lambda\in \R$  the existence of a positive solution to \eqref{refuge-eq-steady} is conditioned by the sign of $\mu_p(\lb{\O}+a_0+\lambda a_1)$.

First let us observe that for $\lambda> \frac{\|k\|_{\infty}-\|a_0\|_{\infty}}{\|a_1\|_{\infty}}$ we have $sup_{x\in \O}(a_0(x)+\lambda a_1(x) -k(x))>0$ and by $(v)$ of Proposition \ref{refuge-prop-mup} we have  $\mu_p(\lb{\O}+a_0+\lambda a_1) \le -sup_{\O}(a_0(x)+\lambda a_1(x) -k(x)) < 0$.  Therefore by  Theorem \ref{refuge-thm-existence}, there exists a positive solution to \eqref{refuge-eq-steady} for all $\lambda>\frac{\|k\|_{\infty}-\|a_0\|_{\infty}}{\|a_1\|_{\infty}}.$ 
Let us consider the following set $\{ \lambda\,|\,\mu_p(\lb{\O}+a_0+\lambda a_1) = 0\}$. 
 When $\{ \lambda\,|\,\mu_p(\lb{\O}+a_0+\lambda a_1) = 0\}\neq \emptyset$, by   monotonicity of $\mu_p$ with respect to $\lambda$ ($(ii)$ of Proposition \ref{refuge-prop-mup} ), we can see that $\{\lambda\,|\,\mu_p(\lb{\O}+a_0+\lambda a_1) = 0\}$ is bounded from above. 
 Therefore we can define $\lambda^*\in [-\infty,+\infty)$ by the following formula
 $$\lambda^* := \sup\{ \lambda\,|\,\mu_p(\lb{\O}+a_0+\lambda a_1) = 0\},$$
 where we set $\lambda^*=-\infty$ when $\{ \lambda\,|\,\mu_p(\lb{\O}+a_0+\lambda a_1) = 0\}=\emptyset$. 
 
 By construction, thanks to Theorem \ref{refuge-thm-existence} for all $\lambda > \lambda^*$ there exists a unique positive continuous solution  to \eqref{refuge-eq-steady} and when $\lambda^*\in \R$, there is no positive solution to \eqref{refuge-eq} for all $\lambda \le \lambda^*$. 
 
 Before proving the trichotomy,  let us look at the asymptotic behaviours with respect to $\lambda$  of the unique solution $u_\lambda$. First let us observe that the map $\lambda \mapsto  u_\lambda$ is monotone non decreasing. Indeed, thanks to the nonnegativity of $a_1$, for any $\lambda\ge \lambda'$, the continuous bounded function $u_{\lambda'}$ is a subsolution of the problem
 \begin{equation}\label{refuge-eq-asymp1}
  \int_{\O} K(x,y)u(y)\, dy - k(x)u(x) + a_0 (x)u + \lambda a_1 (x)u -\beta u^p = 0.
 \end{equation}
 Observe that any large constant M is a super-solution of \eqref{refuge-eq-asymp1}. Therefore by taking $M$ large enough we have $u_{\lambda'} \le M$ and by the monotone iteration scheme we can construct a positive bounded solution of \eqref{refuge-eq-asymp1} which satisfies $u_{\lambda'}\le  u\le M$. We conclude by using the uniqueness of the solution of problem \eqref{refuge-eq-asymp1}. Hence, $u_{\lambda'}\le  u_\lambda \equiv u$.

 The asymptotic behaviour of $u_\lambda$  when $\lambda \to +\infty$ is obtained by establishing a bound from below for the solution $u_\lambda$ when $\lambda \to  +\infty$. More precisely we show that for all $x \in \O_1$ we have for $\lambda$ large enough
 \begin{equation}
 \label{refuge-eq-asymp2} 
 u_\lambda (x) \ge \left(\frac{\lambda a_1 +a_0 -k(x)}{\sup_{\O}\beta}\right)^{\frac{1}{p-1}}.
\end{equation}
 Indeed from \eqref{refuge-eq-steady} using that $u_\lambda$ is non negative we have 
 $$\beta(x) u_\lambda^p (x) \ge [k(x) + a_0 (x) + \lambda a_1 (x)]u_\lambda.$$
  Thus for $x\in \O_1 $ \eqref{refuge-eq-asymp2}  holds for $\lambda$ large enough. From \eqref{refuge-eq-asymp2} we get trivially that for all $x \in \O_1$
  $$   \lim_{\lambda \to +\infty} u_\lambda(x)\ge \lim_{\lambda \to +\infty}\left(\frac{\lambda a_1 +a_0 -k(x)}{\sup_{\O}\beta}\right)^{\frac{1}{p-1}}=+\infty. $$
 
 So for $x \in \O\setminus\O_1$ so that $|B_{\eps_0} (x) \cap \O_1| > 0$ where $\eps_0$ is given by \eqref{refuge-hyp2} we conclude that
 $$ \lim_{\lambda \to +\infty}\int_{B_{\eps_0}(x)\cap \O_1} u_\lambda(y)\, dy=+\infty.$$
 
 Therefore from \eqref{refuge-hyp2}, \eqref{refuge-eq-steady} and $u_\lambda \ge 0$ we deduce that
 $$\left(\beta(x)u_\lambda^{p-1} + k(x) \right)u_\lambda(x) \ge \int_{\O}K(x,y)u_{\lambda}(y)\,dy\ge  c_0\int_{B_{\eps_0}(x)\cap\O_1}u_\lambda(y)\,dy $$
 which leads to
 $$\lim_{\lambda \to +\infty}    u_\lambda \left(\beta(x)u_\lambda^{p-1} + k(x) \right)\ge \lim_{\lambda \to +\infty}c_0\int_{B_{\eps_0}(x)\cap \O_1} u_\lambda(y)\, dy =+\infty \quad  \text{ for all } \quad x \in \bigcup_{z \in \O_1}B_{\eps_0}(z).$$

 The later implies that
$$ \lim_{\lambda \to +\infty}    u_\lambda (x) = +\infty \quad \text{ for all }\quad x \in \bigcup_{z \in \O_1}B_{\eps_0}(z).$$
 By repeating the above argument with $ \bigcup_{z \in \O_1}B_{\eps_0}(z)$ instead $\O_1$, we show that
$$ \lim_{\lambda \to +\infty}    u_\lambda (x) = +\infty \quad \text{ for all }\quad x \in \bigcup_{z \in \O_1}B_{2\eps_0}(z).$$
 By a finite iteration of the above argumentation, we get 
 $$ \lim_{\lambda \to +\infty}    u_\lambda (x) = +\infty \quad \text{ for all }\quad x \in \bar \O.$$

 Let us now deal with the limit of $u_\lambda$ when $\lambda \to \lambda^{*,+}$. First let us assume that $\lambda^{*,+}\in \R$. In this situation by using the positivity of $u_\lambda$ and the monotonicity of $u_\lambda$ with respect to $\lambda$, we deduce that $u_\lambda$ converges pointwise to $u_{\lambda^{*,+}}$ when $\lambda \to \lambda^{*,+}$. Moreover thanks to the Lebesgue dominated convergence Theorem by passing to the limit in \eqref{refuge-eq-steady}, we see that  $u_{\lambda^{*,+}}$ is a non negative solution of \eqref{refuge-eq-steady} with $\lambda = \lambda^{*,+}$. Therefore by Theorem \ref{refuge-thm-existence} we deduce that $u_{\lambda^{*,+}}\equiv 0$ since $\mu_p(\lb{\O}+a_0+\lambda a_1) = 0$. Thus in this case 
  $$ \lim_{\lambda \to \lambda^{*,+}}    u_\lambda (x) = 0 \quad \text{ for all }\quad x \in \bar \O.$$
 Lastly assume that $\lambda^{*,+}=-\infty$. Again by using the positivity of $u_\lambda$ and the monotonicity of $u_\lambda$ with respect to $\lambda$, we deduce that $u_\lambda$ converges pointwise to $u_\infty$ when $\lambda \to -\infty$. Now observe that by the monotonicity of $u_\lambda$, we have for all $\lambda\le 0, u_\lambda \le M_0 := \|u_0\|_{\infty}$ and 
 
 $$ a_1 (x)|\lambda|u_\lambda \le  C_0,$$ where 
 $$C_0 := M_0\left\|\int_{\O}K(\cdot,y)dy + k + a_0\right\|_{\infty}+ \|\beta\|_{\infty} M_0^p.$$ 
  Therefore for $x \in \O_1$ , we deduce that
 $$ 0 \le  u_{\infty} (x) = \lim_{\lambda \to -\infty} u_\lambda (x) \le  \lim_{\lambda \to -\infty}\frac{C_0}{ a_1 (x)|\lambda|} = 0.$$
 By passing to the limit in the equation \eqref{refuge-eq-steady}, thanks to the Lebesgue dominated convergence Theorem we see that $u_\infty$ satisfies the equation below

  \begin{equation}\label{refuge-eq-asymp4}
  \oplb{u}{\O\setminus \O_1} + a_0 (x)u + \lambda a_1 (x)u -\beta u^p = 0.
 \end{equation}

 By Theorem \ref{refuge-thm-existence}, the existence of a positive solution to the above equation is  governed by the sign of $\mu_p(\lb{\O\setminus \O_1} +a_0)$. Therefore when $\mu_p(\lb{\O\setminus \O_1} +a_0) < 0$ there is a unique positive solution whereas for $\mu_p(\lb{\O\setminus \O_1} +a_0) \ge 0$ there is none. In the later case, we deduce that
$$ \lim_{\lambda \to -\infty}    u_\lambda (x) = 0 \quad \text{ for all }\quad x \in \bar \O.$$

Now let us look more closely at the properties of $\lambda^*$ and prove the trichotomy

  \begin{enumerate}
 \item $\lambda^*=-\infty$ when $\mu_p(\lb{\O\setminus \O_1} +a_0) < 0,$
 \item $\lambda^*\in [-\infty,\infty)$ when  $\mu_p(\lb{\O\setminus \O_1} +a_0) = 0,$
 \item $\lambda^*\in \R$ when $\mu_p(\lb{\O\setminus \O_1} +a_0) > 0$.
 \end{enumerate}

 \subsection*{Case 1: $\mu_p(\lb{\O\setminus \O_1} +a_0) < 0$}
 In this situation, observe that  by $(i)$ of Proposition \ref{refuge-prop-mup},  we have for all $\lambda$
 $$0>\mu_p(\lb{\O\setminus \O_1} +a_0) = \mu_p(\lb{\O\setminus \O_1} +a_0+\lambda a_1) \ge \mu_p(\lb{\O}+a_0+\lambda a_1).$$
 Therefore, thanks to Theorem \ref{refuge-thm-existence} there exists a positive non trivial solution to \eqref{refuge-eq-steady} for all $\lambda \in \R$ .  Thus $\lambda^*=-\infty$.
 \subsection*{Case 2:$\mu_p(\lb{\O\setminus \O_1} +a_0) = 0$}
 In this situation, by  monotonicity of $\mu_p$ with respect to $\lambda$ ($(ii)$ of Proposition \ref{refuge-prop-mup} ) and $(i)$ of   Proposition \ref{refuge-prop-mup} either $\mu_p(\lb{\O}+a_0+\lambda a_1) < 0$ for all $\lambda\le 0$ or there exists $\lambda_0\le 0$ so that $\mu_p(\lb{\O}+a_0+\lambda_0 a_1) = 0.$ In the first situation, as above  there exists a positive solution to \eqref{refuge-eq-steady} for any $\lambda$ and $\lambda^*=-\infty$.
 In the other case, $\lambda^*\ge \lambda_0$ and $\lambda^*\in\R$.

  \subsection*{Case 3:$\mu_p(\lb{\O\setminus \O_1} +a_0) >0$}
  In this last situation, we claim that 
  \begin{cla}
 $$\liminf_{\lambda \to -\infty} \mu_p(\lb{\O} +a_0+\lambda a_1) > 0.$$
\end{cla} 
  Assume the claim holds true then this  implies that $\{ \lambda\,|\,\mu_p(\lb{\O}+a_0+\lambda a_1) = 0\}$ is non empty and therefore $\lambda^*\in\R$. Indeed, since  $\mu_p(\lb{\O} +a_0+\lambda a_1) < 0$ for any $\lambda >  \frac{\|k\|_{\infty}-\|a_0\|_{\infty}}{\|a_1\|_{\infty}}$ and by the claim there exists $\bar \lambda$ so that $\mu_p(\lb{\O} +a_0+\bar \lambda a_1) > 0$,  by  continuity of $\mu_p$ with respect to $\lambda$ there exists a $\bar \lambda<\lambda_0 <\frac{\|k\|_{\infty}-\|a_0\|_{\infty}}{\|a_1\|_{\infty}}$ so that $\mu_p (\lb{\O} + a_0 + \lambda_0 a_1 ) = 0.$
  
 \dem{Proof of the Claim:}
 The proof of this claim relies on the construction of an adequate test function. By arguing as in the proof of  Claim  \ref{refuge-cla-super}  we can introduce a regularisation $a_\eps$ of $a_0 -k$ so that the following operator
 $$\oplb{u}{\eps,\O\setminus\O_1}:= \int_{\O\setminus\O_1}K(x,y)u(y)dy + a_\eps (x)u$$
 has a positive continuous principal eigenfunction. By  continuity of $\mu_p (\lb{\eps,\O\setminus\O_1})$ with respect to $a_\eps$ ($(iii)$ of Proposition \ref{refuge-prop-mup}) we can find $\eps$ small so that
 $$\mu_p (\lb{\eps,\O\setminus\O_1} ) -\|a_\eps - a_0 + k\|_{\infty}\ge \frac{\mu_p (\lb{\O\setminus\O_1} + a)}{ 2}.$$
Let $\eps$ be fixed and let  $\O_\delta$ be the set 
 $$\O_\delta := \{x \in \O_1\;|\;  d(x;\partial\O_1) > \delta\}.$$

As in the proof of the Claim \ref{refuge-cla-super},  by  continuity of  $\sup_{\O\setminus \O_\delta}$ and $\mu_p(\lb{\eps,\O\setminus\O_\delta})$  with respect to the domain  we achieve for $\delta$ small enough, say $\delta \le  \delta_0,$ 
 \begin{equation}\label{refuge-eq-liminf1}
 \mu_p (\lb{\eps,\O\setminus\O_1} ) \ge \mu_p (\lb{\eps,\O\setminus\O_\delta} ) \ge \mu_p (\lb{\eps,\O\setminus\O_{\delta_0}} ) \ge \|a_\eps - a_0+k\|_{\infty} +\frac{ \mu_p (\lb{\O\setminus\O_1} + a_0)}{ 8}.
 \end{equation}

By construction $\bar\O\setminus \O_ {\frac{\delta}{2}}$ and $\bar \O_{\delta}$ are two disjoints bounded closed set, so by the Urysohn Lemma there exists a nonnegative continuous function $\eta_1$ such that $0\le \eta_1\le 1, \eta_1 (x) = 1$ in $\bar \O\setminus \O_ {\frac{\delta}{2}} , \eta_1 (x) = 0$ in $\bar\O_{\delta}.$  
  
Let $\psi_1,\psi_2$ be the  following continuous functions
 $$\psi_1 :=\begin{cases} 
 &\Psi_\delta \eta_1 \quad\text{ in }\quad \O\setminus\O_{\delta}\\
 &0 \quad\text{elsewhere},
\end{cases} 
\qquad
 \psi_2 :=\begin{cases} 
&c_0\eta_2\quad \text{ in }\quad \O_1\\
 &0 \quad\text{elsewhere}.
 \end{cases}
 $$
 where $\Psi_\delta$ is the positive continuous eigenfunction associated with $\mu_p (\lb{\eps,\O\setminus\O_\delta} )$ normalized by $\|\Psi_\delta\|_{\infty}= 1$ and $c_0$ is positive constant to be specified later on. 
Consider now the function $\psi:= \sup (\psi_1 ,\psi_2 )$ and let  $\gamma$  be a positive constant to be fixed later on.
 We will prove that for $\gamma,\delta, \lambda$ and $c_0$ well chosen the function $\psi$ is an adequate test function for $\lb{\O}+a_0+\lambda a_1+\gamma$. So let us 
 compute $\oplb{\psi}{\O}  + a\psi +\lambda a_1 \psi +\gamma \psi$. 

On $\O\setminus \O_{\frac{\delta}{2}}$, by construction we have   
\begin{multline*}
 \oplb{\psi}{\O}  + a_0\psi+\lambda a_1\psi+\gamma \psi \le \int_{\O\setminus\O_{\frac{\delta}{2}}}K(x,y)\Psi_\delta(y) \, dy + \| K(\cdot,\cdot)\|_{\infty}\left(|\O_{\frac{\delta}{2}}\setminus\O_\delta| +c_0|\O_\delta|\right)\\
 +a_\eps \Psi_\delta +\|a_\eps +k -a_0\|_{\infty} \Psi_\delta+\lambda a_1\Psi_\delta+\gamma \Psi_\delta.
 \end{multline*}

 Therefore by \eqref{refuge-eq-liminf1} we see that  
 \begin{equation}\label{refuge-eq-liminf2}
  \oplb{\psi}{\O}  + a_0\psi+\lambda a_1\psi+\gamma \psi \le \left(-\frac{\mu_p(\lb{\O\setminus \O_1}+a_0)}{8}+\gamma\right)\Psi_\delta+ \| K(\cdot,\cdot)\|_{\infty}\left(|\O_{\frac{\delta}{2}}\setminus\O_\delta| +c_0|\O_\delta|\right).
  \end{equation} 
 Let $m_0$ be the  constant 
 $$ m_0:=\inf_{\delta \in [0,\delta_0]}\left(\inf_{x\in \O\setminus\O_\delta}\Psi_\delta(x)\right).$$
 We have $m_0>0$ since for all $\delta \in [0,\delta_0]$ $\Psi_\delta$ is positive and continuous in $\bar \O\setminus \O_\delta$. 
 
 Now let us fix $\gamma:=\frac{\mu_p(\lb{\O\setminus\O_1}+a_0)}{16}$ and choose  $\delta$ and $c_0$ such that  
 
\begin{align} 
&|\O_{\frac{\delta}{2}}\setminus\O_{\delta} |\le \frac{m_0\mu_{p}(\lb{\O\setminus\O_1}+a_0)}{64\|K(\cdot,\cdot)\|_{\infty}},\label{refuge-eq-liminf3}\\
&c_0\le \frac{m_0\mu_{p}(\lb{\O\setminus\O_1}+a_0)}{64\|K(\cdot,\cdot)\|_{\infty}|\O_1|}. \label{refuge-eq-liminf4}
 \end{align}
 By combining \eqref{refuge-eq-liminf2}, \eqref{refuge-eq-liminf3}  and \eqref{refuge-eq-liminf4} we see that  
$$\left(-\frac{\mu_p(\lb{\O\setminus \O_1}+a_0)}{8}+\gamma\right)\Psi_\delta+ \| K(\cdot,\cdot)\|_{\infty}\left(|\O_{\frac{\delta}{2}}\setminus\O_\delta| +c_0|\O_\delta|\right) \le 0. $$ 
   
 Therefore  on $\O\setminus\O_{\frac{\delta}{2}}$, we achieve for all $\lambda \le 0$
 \begin{equation}\label{refuge-eq-liminf5}
 \oplb{\psi}{\O}+(a_0+\lambda a_1+\gamma)\psi \le 0. 
 \end{equation}

 Now on $\O_{\frac{\delta}{2}}$  we have by construction 
 $$\oplb{\psi}{\O}+(a_0+\lambda a_1 +\gamma)\psi\le \| K(\cdot,\cdot)\|_{\infty}|\O|+\|a_0\|_{\infty}+\gamma +\lambda c_0\left(\inf_{\O_{\frac{\delta}{2}}}a_1\right). $$
 Since $a_1>0$ in  $\bar \O_{\frac{\delta}{2}}$, by choosing  $\lambda \le -\frac{\| K(\cdot,\cdot)\|_{\infty}|\O|+\|a_0\|_{\infty}+\gamma}{c_0\left(\inf_{\O_{\frac{\delta}{2}}}a_1\right)}$ we get
 \begin{equation} \label{refuge-eq-liminf6}
\oplb{\psi}{\O}+(a_0+\lambda a_1 +\gamma)\psi\le 0 \quad  \text{ in }\quad \O_{\frac{\delta}{2}}.
\end{equation}
 Hence from \eqref{refuge-eq-liminf5} and \eqref{refuge-eq-liminf6} we see that  for $\lambda$ negative enough, the function $(\psi, \gamma)$ is an adequate test function for the operator $\lb{\O} + a_0 + \lambda a_1$. That is:  $\psi$ is a positive continuous function on $\O$ which satisfies 
 $$ \oplb{\psi}{\O}+(a_0+\lambda a_1 +\gamma)\psi\le 0.$$
 So by  definition of $\mu_p(\lb{\O} +a_0+\lambda a_1)$ we deduce that for $\lambda$ negative enough we have
 $\mu_p(\lb{\O} +a_0+\lambda a_1) \ge \gamma  > 0.$ 
 \fdem


\section{The partially controlled problem: The refuge case}\label{refuge-section-refuge}

In this Section, we analyse \eqref{refuge-eq} in the presence of a refuge zone, i.e. when there exists $\o\subset \O$ so that $\beta_{|\o}\equiv 0$.
 In a presence of a refuge zone, the analysis of  \eqref{refuge-eq} is more involved and the characterisation of the existence/non-existence of a positive solution of \eqref{refuge-eq} cannot always be summarised to a single  critical value $\lambda^*$. In this situation, we prove the Theorem \ref{refuge-thm3} that we recall below:
  
\begin{theorem}
 Assume that $K,a_i$ and $\beta$ satisfy \eqref{refuge-hyp1}--\eqref{refuge-hyp2}. Assume further that $\o\neq\emptyset,$  then there exists two quantities $\lambda^*,\lambda^{**} \in [-\infty,+\infty]$ so that  we have the following dichotomy :
\begin{itemize}
\item Either $\lambda^{**}\le \lambda^{*}$ and there exists no positive bounded  solution to  \eqref{refuge-eq}.
\item Or  $\lambda^*< \lambda^{**}$, and  for all $\lambda \in (\lambda^*,\lambda^{**})$ 
 there exists a unique positive bounded continuous solution  to \eqref{refuge-eq}. When $\lambda^*, \lambda^{**}\in \R$, there is no positive bounded solution to \eqref{refuge-eq} for all $\lambda \le \lambda^*$ and for all $\lambda\ge \lambda^{**}$. Moreover, the map $\lambda\to u_\lambda$ is monotone increasing and we have
\begin{itemize}
\item[(i)] 
$$ \lim_{\lambda \to \lambda^{**,-}} \|u_\lambda\|_{\infty,\o}= +\infty, $$
where $\|u_\lambda\|_{\infty,\o}:=\sup_{x\in\o}|u_\lambda(x)|$.
\item[(ii)] If $\mu_p(\lb{\o}+a_0+\lambda^{**}a_1)$ is an eigenvalue in $L^1(\o)$ or $\lambda^{**}=+\infty$ then  $$ \forall\, x \in \bar \O, \lim_{\lambda \to \lambda^{**,-}} u_\lambda (x) = +\infty.$$
\item[(iii)] For all $x \in \bar \O$ we have $ \lim_{\lambda \to \lambda^{*,+}} u_\lambda(x) = u_\infty(x),$
where $u_\infty$ is a  function satisfying on $\O_1 := \{x\in \O| \,a_1 (x) > 0\, u_\infty\equiv 0$ and $u_\infty$ is a nonnegative solution to 
 $$\int_{\O\setminus \O_1} K(x,y)u(y)dy - k(x)u + a_0 (x)u -\beta u^p = 0 \quad \text{in}\quad \O\setminus\O_1.$$
 Furthermore, $u_\infty$ is non trivial when $\mu_p(\lb{\O\setminus \O_1}+a_0)<0$.
\end{itemize}
 \end{itemize}
\end{theorem}

\dem{Proof:}
Thanks to Theorem \ref{refuge-thm-refuge}  the existence of a positive unique bounded  solution to \eqref{refuge-eq} in  presence of a refuge zone  is conditioned by the following inequality 
\begin{equation}\label{refuge-eq-cond1}
 \mu_p(\lb{\o}+a_0+\lambda a_1)>0> \mu_p(\lb{\O}+a_0+\lambda a_1).
 \end{equation}

Let us introduce the following quantities:

\begin{align*}
  \lambda^* := \sup\{ \lambda\,|\,\mu_p(\lb{\O}+a_0+\lambda a_1) \ge 0\},\\
  \lambda^{**}:=\inf\{\lambda\,|\, \mu_p(\lb{\o}+a_0+\lambda a_1) \le 0\}.
  \end{align*}
  
  We can see that the description of the set of positive bounded solutions of \eqref{refuge-eq} is then equivalent to show whether or not we have  $\lambda^*<\lambda^{**}$.  To answer this question,  we analyse separately the three different situations :
 \begin{enumerate}
 \item $\o \subset \O\setminus \O_1$,
 \item $\o \subset \subset \O_1$, 
 \item $\o \not\subset  \O\setminus \O_1$ and $\o \not\subset \O_1$.
 \end{enumerate} 

Let us start with the analysis the first  situation. 
\subsection*{Case 1: $\o \subset \O\setminus \O_1$}
In this situation, since $\o \subset\O\setminus \O_1$,  we have $\mu_p(\lb{\o}+a_0+\lambda a_1)= \mu_p(\lb{\o}+a_0)$. So from \eqref{refuge-eq-cond1}   we see that for all $\lambda$ there is no bounded  solution to \eqref{refuge-eq-steady} when $\mu_p(\lb{\o}+a_0)\le 0$  whereas    the existence of a bounded solution will be conditioned only by 
the sign of $\mu_p(\lb{\O}+a_0+\lambda a_1)$ when $\mu_p(\lb{\o}+a_0)>0$.   In the later case, the analysis of the Section \ref{refuge-section-kpp} can be reproduced, so we get $\lambda^*\in [-\infty,\infty)<\lambda^{**}=+\infty$.

\subsection*{Case 2: $\o \subset \subset\O_1$ }
 
In this situation, since $\o \subset \O_1$, by $(v)$ of Proposition \ref{refuge-prop-mup} we can see that for some positive constant $C$
$$ -\lambda \sup_{\o}a_1 -C     \le \mu_p(\lb{\o}+a_0+\lambda a_1)\le C -\lambda \sup_{\o}a_1.$$
Therefore,  we see that $\lambda^{**}\in \R$ and by definition of $\lambda^{**}$ and $(i)$ of Proposition \ref{refuge-prop-mup} we have 
$$ \mu_p(\lb{\O}+a_0+\lambda^{**} a_1)\le \mu_p(\lb{\o}+a_0+\lambda^{**} a_1)=0.$$

From the above inequality, we get the following dichotomy:
\begin{itemize}
\item Either $\mu_p(\lb{\O}+a_0+\lambda^{**} a_1)<0$ and by $(ii)$ of Proposition \ref{refuge-prop-mup}, we deduce that  $\lambda^*<\lambda^{**}$.

\item Or $\mu_p(\lb{\O}+a_0+\lambda^{**} a_1)=0$ and by definition of $\lambda^*$ we have $\lambda^*\ge \lambda^{**}$ and for all $\lambda$ there is no positive bounded solution to \eqref{refuge-eq}.
\end{itemize}

\subsection*{Case 3:  $\o \not\subset  \O\setminus \O_1$ and $\o \not\subset \O_1$. }
In this situation, since $\o\cap \O_1 \neq \emptyset$, by $(v)$ of Proposition \ref{refuge-prop-mup} we can see that for some positive constant $C$
$$  \mu_p(\lb{\o}+a_0+\lambda a_1)\le C -\lambda \sup_{\o\cap \O_1}a_1.$$
 Therefore $\lambda^{**} \in [-\infty,+\infty)$.
 Now let us observe that in this situation we also have  for all $\lambda$ 
 $$ \mu_p(\lb{\o }+a_0+\lambda a_1)\le \mu_p(\lb{\o\cap (\O\setminus\O_1)}+a_0+\lambda a_1)=\mu_p(\lb{\o\cap (\O\setminus\O_1)}+a_0).$$
 We then have two case to analyse:
 \begin{itemize}
 \item[({\it i})] $\mu_p(\lb{\o\cap (\O\setminus\O_1)}+a_0)\le 0:$\\
    In this situation, from the above inequality, we can already conclude that $\lambda^{**}=-\infty$. Hence in this situation, for all $\lambda$ there is  no positive bounded solution of \eqref{refuge-eq}. 
 \item [({\it ii})]$\mu_p(\lb{\o\cap (\O\setminus\O_1)}+a_0)> 0:$\\
 In this situation, by working as in Section \ref{refuge-section-kpp} we see that there exists $\lambda <<-1$ so that $ \mu_p(\lb{\o }+a_0+\lambda a_1)>0$. Therefore,  $\mu^{**}\in \R$ and we can argue as in the  Case 2. 
 \end{itemize}
 
Let us now look at the asymptotic behaviour of $u_\lambda$ with respect to $\lambda$. The monotone behaviour of $u_\lambda$ and  its limit  as $\lambda \to \lambda^*$ (i.e (iii)) can be obtained by following the arguments  in Section \ref{refuge-section-kpp}, so we drop the proof here and  prove only $(i)$ and $(ii)$ i.e. we analyse the limits of $u_\lambda$ as $\lambda \to \lambda^{**}$.

When $\lambda^{**}=+\infty$, the behaviour of $u_\lambda$ can be obtained by reproducing the arguments of Section \ref{refuge-section-kpp} and we get for all $x\in \bar \O$
$$\lim_{\lambda \to +\infty}u_\lambda(x)=+\infty.$$

 Now, let us assume that $\lambda^{**}\in \R$. 
By definition of $\lambda^{**}$, we must have 
$$\mu_p(\lb{\o}+a_0+\lambda^{**}a_1)=0.$$
As a consequence we also have $\mu_p (\lbs{\o} + a_0+\lambda^{**}a_1)=0$ where $ \lbs{\o}+  a_0+\lambda^{**}a_1$  is defined by
 $$
 \oplbs{\phi}{\o} +  a_0\phi+\lambda^{**}a_1\phi := \int_{\o} K(y,x)\phi(y)dy -k(x)\phi  + a_0\phi+\lambda^{**}a_1\phi.$$

Let us start with the proof of $(i)$.   First assume that $\stackrel{\circ}{\o}\neq \emptyset$, then
by Theorem \ref{refuge-thm-mu} there exists a positive measure  $d\mu^*$ associated with $\mu_p(\lbs{\o}+a_0+\lambda^{**}a_1)$. Moreover, $\phi_p^*(x)dx$ the regular part of $d\mu^*$ satisfies
$\inf_{\o}\phi_p^*>0$.
By integrating the equation \eqref{refuge-eq} over $\o$ with respect to the measure $d\mu^*$, we get  for any $\lambda^*<\lambda<\lambda^{**}$
$$ \int_{\o}\left(\int_{\O}K(x,y)u_\lambda(y)\,dy\right)d\mu^*+\int_{\o}(-k(x)+a_0+\lambda a_1)u_\lambda \,d\mu^*=0.$$

By definition of $\mu_p(\lbs{\o}+a_0+\lambda^{**})$, it follows  from the above equality  
$$ \int_{\o}\left(\int_{\O\setminus \o}K(x,y)u_\lambda(y)\,dy\right)d\mu^*=(\lambda^{**}-\lambda)\int_{\o} a_1u_\lambda \,d\mu^*.$$
Therefore, by using the monotonicity of the map $u_\lambda$ we have for $\lambda_0<\lambda<\lambda^{**}$
$$\frac{1}{\lambda^{**}-\lambda} \int_{\o}\left(\int_{\O\setminus \o}K(x,y)u_{\lambda_0}(y)\,dy\right)d\mu^*\le\|a_1\|_{\infty} \|u_\lambda\|_{\infty,\o}\int_{\o}d\mu^*, $$
which enforces 
$$ \lim_{\lambda \to \lambda^{**}}\|u_{\lambda}\|_{\infty,\o}=+\infty.$$

Assume now that $\stackrel{\circ}{\o}=\emptyset$. In this situation, we have $\mu_p(\lb{\o}+a_0+\lambda^{**}a_1)=-\sup_{\o}(-k+a_0+\lambda^{**}a_1)$. Since $\o$ is a compact set there exists $x_0 \in \o$ so that $-k(x_0) +a_0(x_0)+\lambda^{**}a_1(x_0)=0$. 
We can check that $x_0\in \O_1$, otherwise we have $\sup_{\o\cap(\O\setminus\O_1)}(-k+a_0)=\sup_{\o\cap(\O\setminus\O_1)}(-k+a_0+\lambda^{**} a_1)= 0$ and $\mu_{p}(\lb{\o\cap (\O\setminus\O_1)}+a_0)=0$. The latter equality leads to the contradiction $-\infty<\lambda^{**}=-\infty$, since for all $\lambda$ we have 
$$\mu_{p}(\lb{\o}+a_0+\lambda a_1)\le \mu_{p}(\lb{\o\cap (\O\setminus\O_1)}+a_0)  =0.$$
Now at $x_0$, we have 

$$\int_{\O}K(x_0,y)u_\lambda(y)\,dy=(\lambda^{**}- \lambda)a_1(x_0)u_\lambda(x_0). $$

By using that $u_\lambda$ is monotone with respect to $\lambda$ we get for all $\lambda_0 \le  \lambda <\lambda^{**}$ 
$$\frac{1}{(\lambda^{**}-\lambda)a_1(x_0)}\int_{\O}K(x_0,y)u_{\lambda_0}(y)\,dy=u_\lambda(x_0), $$
which implies 
$$ \lim_{\lambda \to \lambda^{**}}u_{\lambda}(x_0)=+\infty.$$

Let us now prove $(ii)$. When $\mu_p(\lb{\o}+a_0+\lambda^{**}a_1)$ is associated with a positive $L^1(\o)$ eigenfunction  we claim that 

\begin{cla}
$$\lim_{\lambda \to \lambda^{**}}\int_{\O}u_\lambda(x)\, dx=+\infty.$$
\end{cla}
Assume for the moment that the claim holds then we get $(ii)$  by arguing as follows.
Since $\bar \O$ is compact, in view of the claim there exists $\bar x \in \bar \O$ so that 
 $$ \lim_{\lambda \to \lambda^{**}}\int_{B(\bar x,\frac{\eps_0}{4})\cap \O}u_\lambda(x)\,dx =+\infty.$$
From the equation \eqref{refuge-eq} and by the assumption \eqref{refuge-hyp2} we always have 
$$c_0 \int_{B(x,\eps_0)\cap \O}u_\lambda (x)\, dx\le (k(x)-a_0(x)-\lambda a_1(x))u_\lambda (x) +\beta(x)u^p_\lambda(x).$$
Therefore for all $x \in B(\bar x,\frac{\eps_0}{2})$, we have $B(\bar x ,\frac{\eps_0}{4})\subset B(x,\eps_0)$ which combined with the above inequalities implies that 
$$\lim_{\lambda \to \lambda^{**}}(k(x)-a_0(x)-\lambda a_1(x))u_\lambda (x) +\beta(x)u^p_\lambda(x)=+\infty.$$ 
 Thus  $$ \lim_{\lambda \to \lambda^{**}}u_\lambda(x)=+\infty \quad \text{ for all }\quad x\in B(\bar x,\frac{\eps_0}{2})\cap \bar \O.$$
  In  the above arguments by replacing $\bar x$ by  any $x\in B(\bar x,\frac{\eps_0}{4})\cap \O$, we achieve 
  $$ \lim_{\lambda \to \lambda^{**}}u_\lambda(x)=+\infty \quad \text{ for all }\quad x\in \bigcup_{x\in B(\bar x,\frac{\eps_0}{4})}B(x,\frac{\eps_0}{2}) \cap \O.$$
  Since $\O$ is compact we achieve $\lim_{\lambda \to \lambda^{**}}u_\lambda(x)=+\infty $ for all $x\in \O$ after a finite iteration of this argument.   
\fdem

\dem{Proof of the Claim}

Assume by contradiction that $ \sup_{\lambda}\int_{\O}u_\lambda(x)\, dx < +\infty.$
Since $u_\lambda$ is monotone, by Lebesgue  monotone convergence Theorem  we have $u_\lambda \to \bar u$ in $L^1(\O)$ when $\lambda \to\lambda^{**}$ and 
$\bar u> u_{\lambda_0}>0$ satisfies the equation

$$\oplb{\bar u}{\O}(x)+(a_0+\lambda^{**}a_1)\bar u(x) -\beta(x)\bar u^p(x)=0 \quad \text{ for almost every $x$ in } \O.$$
Therefore   we have
\begin{equation}
\oplb{\bar u}{\O}(x)+(a_0+\lambda^{**}a_1)\bar u(x)=0 \quad \text{ for almost every $x$ in }  \o.\label{refuge-eq-cla-asb}
\end{equation}

By assumption there exists a positive $L^1$ eigenfunction $\phi_p$ associated with $\mu_p(\lb{\o}+a_0+\lambda^{**}a_1)$.
Moreover the positive function $\frac{1}{k(x)-a_0(x)-\lambda^{**}a_1(x)}\in L^1(\o)$ and the compact operator $\k$:
$$\begin{array}{ccl}
 C(\o)&\to& C(\o)   \\
      v        &\mapsto &\opk{v}:=\int_{\o}K(x,y)\frac{v(y) dy}{k(y)-a_0(y)-\lambda^{**}a_1(y)}.
     \end{array}$$
is well defined.
By construction, we can check that $\mu_p(\k)=-1$. Indeed, let $v_p:=(k(x)-a_0(x)-\lambda^{**}a_1(x))\phi_p$ then we can see that    $v_p$ is positive and continuous, since by assumption we have  
$$\oplb{\phi_p}{\o}-v_p=0.$$

Moreover, $v_p$ satisfies  $\opk{v_p}=v_p$. Thus by the Krein-Rutman theory, we have $\mu_p(\k)=-1$ and $\psi_1:=v_p$ where $\psi_1$ is the principal positive continuous eigenfunction associated with $\mu_p(\k)$.

Let us now consider $\k^*$ the following compact operator 
$$\begin{array}{ccl}
 L^1(\o)&\to& L^1(\o)   \\
      v        &\mapsto &\opks{v}:=\frac{1}{k(x)-a_0(x)-\lambda^{**}a_1(x)}\int_{\o}K(y,x)v(y) dy.
     \end{array}$$

By the Krein-Rutman Theory there exists  an eigenvalue $\nu_1$ associated with a positive $L^1(\o)$ function $\phi^*$. Furthermore we can check that $\nu_1=-1$. Indeed, since $\phi^*$ is associated with    $\nu_1$  we have 
$$\opks{\phi_p^*}=-\nu_1\phi^*_p.$$
By multiplying the above equation by $v_p$ and then integrating over $\o$ it follows that
\begin{align*} 
\nu_1\int_{o}v_p(x)\phi^*(x)dx&=-\int_{\o}\frac{v_p(x)}{k(x)-a_0(x)-\lambda^{**}a_1(x)}\left(\int_{\o}K(y,x)\phi^*(y) dy\right)dx\\
&=-\int_{\o}v_p(y)\phi^*(y)dy,
\end{align*}

which implies that  $\nu_1=-1$.
  
 Let $\bar v$ be the $L^1(\o)$ function $\bar v:=(k(x)-a_0(x)-\lambda^{**})\bar u$ then we get  by \eqref{refuge-eq-cla-asb}
 \begin{equation}\label{refuge-eq-cla-asb2}
 \opk{\bar v} -\bar v= -\int_{\O\setminus \o}K(x,y)\bar u(y)\,dy \quad \text{ for almost every $x$ in  } \o. 
 \end{equation}
 Since $\opk{\bar v}$ is continuous, we deduce from \eqref{refuge-eq-cla-asb2} that $\bar v$ is continuous in $\o$. 
 So by multiplying \eqref{refuge-eq-cla-asb}   by $\phi^*$ and then integrating   over $\o$ we get 
$$\int_{\o}\phi_1^*(x)\left(\int_{\O\setminus \o}K(x,y)\bar u(y)\,dy\right) dx =\int_{\o}\phi_1^*(x)[\opk{\bar v}-\bar v]dx$$
Since $\nu_1=-1$ and $\bar u>0$ we end up with the contradiction 
$$0<c_0\le \int_{\o}\phi_1^*(x)\int_{\O\setminus \o}K(x,y)\bar u(y)\,dy dx = 0.$$

\fdem

\bibliographystyle{amsplain}
\bibliography{refuge.bib}

\end{document}